\documentclass[12pt]{amsart}
\usepackage{amssymb}

\usepackage{graphicx}

\usepackage[]{amsfonts}

\def\underset#1#2{\mathrel{\mathop{\kern0pt #2}\limits_{#1}}}

\def\overset#1#2{\mathrel{\mathop{\kern0pt #2}\limits^{#1}}}

\def\couleur(#1 #2 #3)
	{
\special{ps::[end]}}

\def\sqr#1#2{{\vcenter{\vbox{\hrule height.#2pt
             \hbox{\vrule width.#2pt height#1pt \kern#1pt
             \vrule width.#2pt}
             \hrule height.#2pt}}}}

\def\st{\mathinner{\mkern1mu\raise1pt\hbox{.}				
		   \mkern1mu\raise4pt\hbox{.}
		   \mkern1mu\raise1pt\hbox{.}
		 }
         }

\def\bx#1{\setbox1=\hbox{\kern3pt{#1}\kern3pt}				
 \dimen1=\ht1 \advance\dimen1 by 3pt \dimen2=\dp1 \advance\dimen2 by 3pt
 \setbox1=\hbox{\vrule height\dimen1 depth\dimen2\box1\vrule}%
 \setbox1=\vbox{\hrule\box1\hrule}%
 \advance\dimen1 by .4pt \ht1=\dimen1
 \advance\dimen2 by .4pt \dp1=\dimen2 \box1\relax}

\def\k#1{\kern#1em}
\def\vci{\vrule  width.02em height1.47ex depth-.0ex}				
\def\11{{\rm\k{.2}\vci\k{-.37}1}}

\begin{document}
\title{Sound traveling-waves in wind instruments as solutions to non linear 
non homogeneous gas dynamics equations}
\author{Alain-Yves LeRoux}
\address{UNIVERSITE BORDEAUX1, Institut Math{\'e}matiques de Bordeaux, 
UMR 5251,351,Cours de la Lib{\'e}ration, 33405, Talence Cedex}
\email{Alain-Yves.Leroux@math.u-bordeaux1.fr\\}
\maketitle
\begin{abstract} {
The sound propagation is usually described by a linear homogeneous wave 
equation, though the air flow in a duct is described by the gas dynamics 
equations, using a variable cross section, which corresponds to a non 
linear non homogneous system. The aim of this paper is to exhibit a common 
periodic solution to both models, with several free parameters such as 
frequency or amplitude, able to represent any sound. By taking in account 
a friction term linked to the material (wood or brass for instance) of 
the duct, it is possible to build an analytic such solution when the 
cross section fullfills some condition which corresponds exactly to the 
general shape of the wind instruments. The conclusion is that in the 
wind intruments, the shape yields linearity.  \ \par
}\end{abstract}
\section{The gas dynamics models with source terms{\hskip 1.8em}}
\setcounter{equation}{0}The duct is supposed to have a cylindrical symetry 
along its axis which can be a straight line (as for a clarinet) or a 
curve line (as for a horn) and the abscissa along this axis is denoted 
by $x>0,$ with $x=0$ corresponding to the narrowest side of the duct. 
The variable cross section is denoted by $a(x)>0,$ the air density and 
velocity are $\rho $ and $u$, both depending on the position $x$ and 
time $t.$ Since people blows at $x=0,$ the velocity $u$ is expected to 
be positive. The mass conservation is described by
\begin{displaymath} 
a(x)\ \rho _{t}\ +\ \left({a(x)\ \rho \ u}\right) _{x}\ =\ 0\ ,\end{displaymath} 
which becomes the transport equation 
\begin{equation} 
q_{t}\ +\ m_{x}\ =\ 0\ ,\label{GasDAcoustics0}
\end{equation} by using the transported quantity $q=a(x)\rho $ and the 
flux $m=a(x)\rho u.$  This flux is ruled by the general dynamic equation
\begin{equation} 
m_{t}\ +\ 2\ u\ m_{x}\ +\ \displaystyle \left({\displaystyle c^{2}-u^{2}}\right) 
\ q_{x}\ +\ S(q,m)\ =\ 0\ ,\ \label{GasDAcoustics1}
\end{equation} where $c=c(q)\ $is the sound velocity and $S$ is a source 
term to be detailed later. The system ~(\ref{GasDAcoustics0}),~(\ref{GasDAcoustics1}) 
is an inhomogeneous form of the isentropic Euler equations, with the 
usual wave velocities $u-c$ and $u+c.$ \ \par
{\hskip 1.5em}A traveling-wave is expected to have the form $\ q(x,t)=q(x-At),\ m=m(x-At)\ ,\ $with 
some constant $A>0$ corresponding to a velocity, so that both $q$ and 
$m$ are solution to the single wave equation:$\ \ \ q_{t}\ +\ A\ q_{x}\ =\ 0\ ,\ m_{t}\ +\ A\ m_{x}\ =\ 0\ .$ 
  \ \par
{\hskip 1.8em}Following the notion of {\it{source waves}}, such that 
Roll waves in hydraulics (see [1] or [3]), such solutions are exhibited 
by writing $m=m(q)$ in ~(\ref{GasDAcoustics0}) and ~(\ref{GasDAcoustics1}). 
One gets 
\begin{displaymath} 
\left({\left({m'(q)-u}\right) ^{2}-c^{2}}\right) q_{x}\ =\ S(q,m(q))\ 
.\end{displaymath} Next, by introducing a function $\psi (q)$ such that 

\begin{displaymath} 
\psi '(q)\ =\ \frac{\left({m'(q)-u}\right) ^{2}\ -\ c^{2}}{S(q,m(q))}\ 
,\end{displaymath} we obtain the two equations $\psi '(q)\ q_{x}\ =\ 1$ 
and $\ \ \ \psi '(q)\ q_{t}\ =\ -m'(q)\ \psi '(q)\ q_{x}\ =\ -m'(q)\ .$\ 
\par
By integrating the first equation with respect to $x$ we get$\ \psi (q)\ =\ x\ -\ K(t)\ ,$ 
where $K(t)$ does not depend on $x$. Next, a derivation with respect 
to $t$ , using $\psi '(q)\ q_{t}=-m'(q)\ ,$ leads to $m'(q)=K'(t)\ .$ 
Another derivation, with respect to $x$ this time gives $m''(q)\ q_{x}=0\ .$ 
Since the expected wave is not flat, we select $m''(q)=0,\ $which leads 
to the two expressions $m(q)\ =\ A\ q\ -\ B\ \ ,\ \ K'(t)\ =\ A\ ,$ with 
some constants $A$ and $B\ ,$ and the relation 
\begin{equation} 
\psi (q)\ =\ x-At\ ,\label{GasDAcoustics2}
\end{equation}  with the integration constant put inside the expression 
of $\psi (q).$ From ~(\ref{GasDAcoustics2}) it is clear that $q$ and 
then $m=Aq-B$ are both functions of $x-At$ , and therefore solutions 
of a linear wave equation with the constant velocity $A$ (traveling-waves). 
  The expression of $\psi '(q)$ gets simpler by noticing that $u=\frac{m}{q}\ =\ A-\frac{B}{q}\ $so 
that $m'(q)-u=\frac{B}{q}\ ,$ which leads to 
\begin{equation} 
\psi '(q)\ =\ \frac{\displaystyle B^{2}\ -\ q^{2}c^{2}}{\displaystyle 
q^{2}\ S(q,Aq-B)}\ .\label{GasDAcoustics3}
\end{equation} It is important to notice that the hypothesis $S(q,m)\not\equiv 0$ 
is fondamental for this result, which does not apply for homogeneous 
systems. Another fondamental hypothesis is that $\psi '(q)$ is a function 
of $q$ {\it{only}}. The result does not apply to the case of a source 
term depending on $q,\ m$ and{\it{ }}$x$ for instance or to a sound speed 
depending on $q$ and $x$. In the next sections we look for state laws 
and duct shapes for which these hypotheses are verified, and conclude 
that it is genuinely the case for the usual shapes of wind instruments, 
as flutes, clarinets or horns.\ \par
\section{Application to an air flow in a duct{\hskip 1.8em}}
\setcounter{equation}{0}The cross section $a(x)$ is supposed to be increasing 
with $x$ ($a'(x)>0$). We look for expressions of $c(q)$ and $S(q,m)$ 
compatible with the previous hypotheses, that is independant of $x$. 
The expression of $c(q)$ comes from the pressure $P$ using $c^{2}=\frac{\gamma P}{\rho }\ $, 
where $\gamma =1.4$ is the adiabatic constant, and the pressure is linked 
to the density $\rho $ and the temperature $T$ (in {\char'27}Kelvin) 
by a state law. A well known state law is given by the Boyle-Mariotte 
law $P=K_{0}\rho T,\ $with $K_{0}=287.06$ in M.K.S. units, with $P$ expressed 
in $Pascals$, which provides the following sound velocity
\begin{displaymath} 
c\ =\ \sqrt{\gamma \ K_{0}\ T}\ .\end{displaymath}  Since the usual temperatures 
are greater than 270{\char'27}K, the variations of the profile of $c$ 
are negligible, and we can fix the temperature to a reference value $T_{0}$ 
and take $c=c_{0}\equiv \sqrt{\gamma \ K_{0}\ T_{0}}\ .$ This corresponds 
to the "{\it{isothermal}}" hypothesis. \ \par
{\hskip 1.8em}Another state law is given by the "{\it{isentropic}}" case 
$P=K\ \rho ^{\gamma }$, with $K=69259.5\ $in M.K.S. units. The expression 
of the sound speed reads $c=\sqrt{\gamma K}\ \rho ^{\frac{\gamma -1}{2}}\ ,$ 
 which is not a function of $q$ only.  Since $c^{2}=P'(\rho ),$ we can 
derive the expression of the source term at rest, that is when the velocity 
is zero and the pressure is constant: $u=0,\ \frac{\partial P}{\partial x}=0\ .$ 
The dynamical equation ~(\ref{GasDAcoustics1}) reduces to $c^{2}q_{x}+S(q,0)=0\ .$ 
Since $q_{x}=a(x)\ \rho _{x}\ +\ a'(x)\ \rho \ ,$ we obtain $c^{2}\ a'(x)\ \rho \ +\ S(q,0)\ =\ 0\ ,$ 
which leads to two possible of the source term:
$$ S(q,0)=
\begin{cases} 
{-\ c_{0}^{2}\ \frac{a'(x)}{a(x)}\ q}&{\ 
in\ the\ isothermal\ case\ ,}\cr 
{-\gamma \ K\ \ \frac{a'(x)}{a(x)^{\gamma}}\ q^{\gamma }}&{\ in\ the\ isentropic\ case\ .}\cr  
\end{cases} 
$$
We introduce a friction term, of the Strickler type, to select the following 
source term
\begin{displaymath} 
S(q,m)\ =\ S(q,0)\ +\ k\ \left\vert{u}\right\vert u\ ,\end{displaymath} 
 where $k$ is the Strickler friction (dimensionless) coefficient.\ \par
{\hskip 1.8em}The retained model is a combination between the isothermal 
and the isentropic cases. The sound velocity is a constant $c=c_{0}\equiv \sqrt{\gamma K_{0}T_{0}},\ $with 
a fixed temperature $T_{0}.$ For instance with $T_{0}=300°K$ we get $c_{0}=347.225\ ms^{-1}.$ 
\ \par
{\hskip 1.5em}The selected source term is $S(q,m)=k\left\vert{u}\right\vert u-\gamma Kq^{\gamma }\ \frac{a'(x)}{a(x)^{\gamma 
}}\ ,$which becomes independent on $x$ if we impose a profile such that
\begin{equation} 
\frac{\displaystyle a'(x)}{\displaystyle a(x)^{\gamma }}\ =\ Constant\ 
=\ k\ \frac{\displaystyle D^{2}}{\displaystyle \gamma K}\ ,\label{GasDAcoustics4}
\end{equation} where $D$ is a constant.  Expression ~(\ref{GasDAcoustics4}) 
is a differential equation whose solutions are 
\begin{displaymath} 
a(x)\ =\ \left({\frac{\gamma \ K}{\left({\gamma -1}\right) \ k\ D^{2}\ 
\left({x_{0}-x}\right) }}\right) ^{\frac{1}{\gamma -1}}\end{displaymath} 
 where $x_{0}$ is a constant to be taken larger than $L$ the length of 
the wind intrument. The radius of the duct is given by $r(x)=\sqrt{\frac{a(x)}{\pi }}\ ,$ 
and we recognize the usual shape of the wind instrument as seen on the 
figure below, where $kD^{2}=75\ \gamma K$ . The value of $x_{0}$ corresponds 
to an infinite radius, or to a maximal length for a wind instrument. 
By the way, using the isothermal case leads to another expression of 
the cross section which reads $a(x)=a_{0}\ exp\left({Cx}\right) $ where 
$C$ and $a_{0}$ are positive constants, and allows an infinite length.\ 
\par

\begin{figure}[h]
\begin{center}
\rotatebox{0}{\resizebox{10cm}{!}{\includegraphics{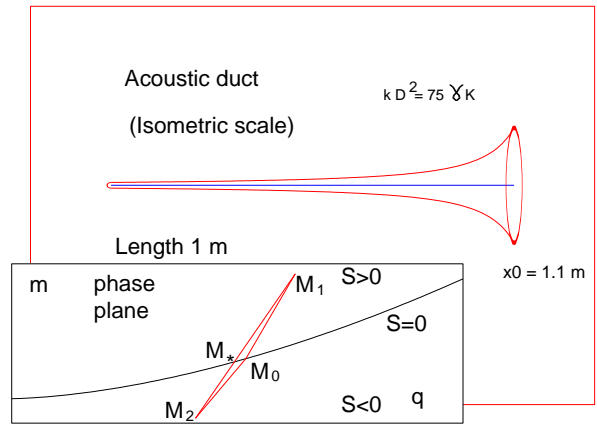}}}
\caption{The shape and the phase plane}
\label{aq0}
\end{center}
\end{figure}
\ \par
{\hskip 1.8em}The dynamic equation ~(\ref{GasDAcoustics1}) becomes
\begin{equation} 
m_{t}+2um_{x}+\displaystyle \left({\displaystyle c_{0}^{2}-u^{2}}\right) 
q_{x}\ +\ k\displaystyle \left({\displaystyle u^{2}-D^{2}q^{\gamma }}\right) 
\ =\ 0\ ,\label{GasDAcoustics5}
\end{equation}  where the source term is independent on $x$ which allows 
waves of the form $m=Aq-B.$ \ \par
\section{The wave profile{\hskip 1.8em}}
\setcounter{equation}{0}The updated form of $\psi '(q)$ reads 
\begin{equation} 
\psi '(q)\ =\ \frac{\displaystyle B^{2}-c_{0}^{2}\ q^{2}}{\displaystyle 
k\displaystyle \left({\displaystyle \displaystyle \left({\displaystyle 
Aq-B}\right) ^{2}\ -\ D^{2}\ q^{\displaystyle \gamma +2}}\right) }\ \label{GasDAcoustics6}
\end{equation}  and we consider a so-called {\it{Reference state}} $M_{*}=\left({q_{*},m_{*}}\right) \ ,$ 
with $q_{*}>0,\ m_{*}>0$ and set $u_{*}=\frac{m_{*}}{q_{*}}\ $and take 
$A=u_{*}+c_{0}$ and $B=q_{*}\ c_{0}$. The value $q_{*}$ is a possible 
realistic value of$q$ which is a root of the numerator in ~(\ref{GasDAcoustics6}). 
Since $\psi '(q)q_{x}=1$, either $q_{x}$ becomes infinite when $q=q_{*},$ 
or $q_{*}$ is also a root of the denaminator in ~(\ref{GasDAcoustics6}) 
and in this case $u_{*}=D\ q_{*}^{\frac{\gamma }{2}}\ ,$ and $M_{*}$ 
belongs to the set $S_{0}=\{\left({q,m}\right) \ \mid \ S(q,m)=0\ \}\ $in 
the phase plane $(q,m)$. The line $m=\left({u_{*}+c_{0}}\right) q-q_{*}c_{0}$ 
cuts $S_{0}$ in $M_{*}$, coming from the set $\{S<0\}$ for $q<q_{*}$ 
and going towards $\{S>0\}$ for $q>q_{*}.$ Since 
\begin{displaymath} 
\psi '(q)\ =\ \frac{\left({q_{*}^{2}-q^{2}}\right) \ c_{0}^{2}}{q^{2}\ 
S\left({q,\left({u_{*}+c_{0}}\right) q-q_{*}c_{0}}\right) }\ ,\end{displaymath} 
we see that for $\psi '(q)<0$ for $q<q_{*},\ S<0$ and for $q>q_{*},\ S>0$. 
Since $\psi '(q)q_{x}=1,$ we get that $q_{x}<0,\ $which means that the 
wave profile corresponding to the line $m=\left({u_{*}+c_{0}}\right) q-q_{*}c_{0}$ 
is decreasing. \ \par
{\hskip 1.8em}Now let us consider in $\{S>0\}$ a point $M_{1}=\left({q_{1},m_{1}}\right) $ 
on the line $m=\left({u_{*}+c_{0}}\right) q-q_{*}c_{0}$ ,and then the 
line whose reference point is $M_{1},$ that is the line of equation $m=\left({u_{1}+c_{0}}\right) q-q_{1}c_{0}$. 
Along this line, we have 
\begin{displaymath} 
\psi '(q)\ =\ \frac{\left({q_{1}^{2}-q^{2}}\right) \ c_{0}^{2}}{q^{2}\ 
S\left({q,\left({u_{1}+c_{0}}\right) q-q_{1}c_{0}}\right) }\ ,\end{displaymath} 
which is positive for $q<q_{1}$. This means that the corresponding profile 
is increasing. By the same way, let us consider in $\{S<0\}$ a point 
$M_{2}=\left({q_{2},m_{2}}\right) $ on the line $m=\left({u_{*}+c_{0}}\right) q-q_{*}c_{0}$ 
,and then the line whose reference point is $M_{2},$ of equation $m=\left({u_{2}+c_{0}}\right) q-q_{2}c_{0}$. 
Along this line, we have 
\begin{displaymath} 
\psi '(q)\ =\ \frac{\left({q_{2}^{2}-q^{2}}\right) \ c_{0}^{2}}{q^{2}\ 
S\left({q,\left({u_{2}+c_{0}}\right) q-q_{2}c_{0}}\right) }\ ,\end{displaymath} 
which is positive for $q>q_{2}$. This means that the corresponding profile 
is increasing too. This is reprensented in a window on both figure, above 
and below. Up to now, the two points $M_{1}$ and $M_{2}$ have been chosen 
without any restriction, and it is possible to get the two lines  $m=\left({u_{1}+c_{0}}\right) q-q_{1}c_{0}$ 
and  $m=\left({u_{2}+c_{0}}\right) q-q_{2}c_{0}$ cutting $S_{0}$ at the 
same point $M_{0}\in S_{0}.$ This can be done by choicing first $M_{0}\in S_{0}\ $with 
$q_{0}>q_{*}$, then searching a point $M=(q,(u_{*}+c_{0})q-q_{*}c_{0})$ 
such that $M_{0}$ belongs to the lines of slope $u_{j}-c_{0}$ passing 
through $M_{j},\ j=1,2.$ This corresponds to the equation
\begin{equation} 
q+\frac{\displaystyle q_{*}q_{0}}{\displaystyle q}\ =\ \displaystyle \left({\displaystyle 
u_{*}+2c_{0}-Dq_{0}^{\displaystyle \frac{\displaystyle \gamma }{\displaystyle 
2}}}\right) \frac{\displaystyle q_{0}}{\displaystyle c_{0}}\ ,\label{GasDAcoustics7}
\end{equation} which have two solutions which are actually $q_{1}$ and 
$q_{2}$ with $q_{2}<q_{*}<q_{0}<q_{1}$ as shown in the figure below.\ 
\par

\begin{figure}[h]
\begin{center}
\rotatebox{0}{\resizebox{9cm}{!}{\includegraphics{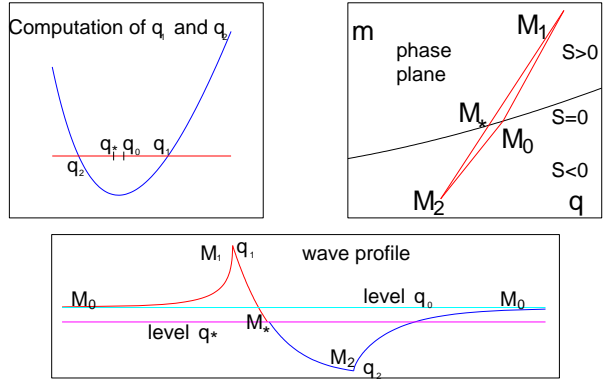}}}
\caption{The profile of the wave}
\label{aq1}
\end{center}
\end{figure}
\ \par
{\hskip 1.8em}The profiles are computed by integrating $\psi '(q)\ q_{x}=1$ 
for each part of the profile, using a numerical integration method, with 

\begin{displaymath} 
\psi '(q_{*})\ =\ \frac{2q_{*}\ c_{0}^{2}}{k\left({2\left({u_{*}+c_{0}}\right) 
u_{*}q_{*}^{2}-D^{2}\left({\gamma +2}\right) q_{*}^{\gamma +1}}\right) 
}\ \end{displaymath} for $q=q_{*}$. We first compute $M_{1}$ and $M_{2}$ 
from ~(\ref{GasDAcoustics7}) and starting from $M_{1}$ we compute the 
two branches $M_{1}M_{0}$ and $M_{1}M_{*}$, next starting from $M_{2}$ 
we compute the two branches $M_{2}M_{0}$ and $M_{2}M_{*}$. \ \par
\section{The wave propagation{\hskip 1.8em}}
\setcounter{equation}{0}The profile is made of a sequence of 3 travelling 
waves: a rear part from the state $M_{0}$ to $M_{1}$ with the constant 
velocity $A_{1}=u_{1}+c_{0}=u_{*}+2c_{0}-c_{0}\frac{q_{*}}{q_{1}}\ ,$ 
a central part from the state $M_{1}$ to $M_{*}$ and $M_{2}$ with the 
constant velocity $A_{*}=u_{*}+c_{0}$ and a front part from the state 
$M_{2}$ back to $M_{0}$ with the velocity $A_{2}=u_{2}+c_{0}=u_{*}+2c_{0}-c_{0}\frac{q_{*}}{q_{2}}\ .$ 
Since $q_{2}<q_{*}<q_{0}<q_{1},$ we have $A_{1}>A_{*}>A_{2}$ . When $q_{*}$ 
and $q_{0}$ are near, which is the case in reality, these velocities 
are also very near and we can consider that the whole wave travels at 
the velocity $u_{*}+c_{0}$. The figure above corresponds to a choice 
of the parameters such that $q_{*}$ and $q_{0}$ are not too near and 
then the three lines in the phase plane do not seem superimposed for 
the reader.\ \par
{\hskip 1.8em}The effect of the different velocities together with the 
mass conservation principle is the emergence of two small shock waves, 
one on the top of the rear wave, in $M_{1}$ and another at the lowest 
part of the front wave in $M_{2}$. The amplitude of these two shocks 
increases during the propagation but stay imperceptible since the travel 
time in the duct is very short, of the order of $10^{-3}\ s$. However, 
such mass effects can be observed in the real world, with similar equations, 
 sometimes at a very different scale in space and time (see for instance 
[4], on the modelling of the Rogue Waves using a small difference of 
the velocities of two long waves on the ocean).  \ \par
{\hskip 1.8em}The sound production is produced by the concatenation of 
such waves with the same parameters, involving a constant wavelength 
and therefore a constant frequency.  \ \par
\section{Some remarks and conclusion{\hskip 1.8em}}
\setcounter{equation}{0}The wave profile presents a decreasing part, 
from $M_{1}$ to $M_{2}$, which never degenerates into a shock as often 
expected in hydrodynamics (as for the Roll waves, see [1] or [3]). This 
comes obviously from the source term. However, the occurrence of a shock 
may be analysed. Since the sound speed is a constant, the pressure law 
corresponds to $P(q)=c_{0}^{2}q$ and the shock condition derived from 
the well known Rankine-Hugoniot relations for a shock wave linking $M_{1}$ 
to $M_{2}$ read
\begin{displaymath} 
\frac{u_{2}-u_{1}}{q_{2}-q_{1}}\ =\ \sqrt{\frac{P(q_{2})-P(q_{1})}{q_{2}-q_{2}}}\ 
=\ \frac{c_{0}}{\sqrt{q_{1}q_{2}}}\ .\end{displaymath} We shall compute 
the rate $\frac{u_{2}-u_{1}}{q_{2}-q_{1}}\ $by two other ways. First 
we write that $M_{0}$ belongs to the lines issued from $M_{1}$ and $M_{2},\ $with 
the respective slopes $u_{1}+c_{0}$ and $u_{2}+c_{0}$, that is 
\begin{displaymath} 
\left({m_{0}\ =\ }\right) \ \left({u_{1}+c_{0}}\right) q_{0}-c_{0}q_{1}\ 
=\ \left({u_{2}+c_{0}}\right) q_{0}-c_{0}q_{2}\ ,\end{displaymath}  which 
gives
\begin{displaymath} 
\frac{u_{2}-u_{1}}{q_{2}-q_{1}}\ =\ \frac{c_{0}}{q_{0}}\ .\end{displaymath} 
 Next we write that $M_{*},\ M_{1}$ ans $M_{2}$ are aligned or 
\begin{displaymath} 
u_{1}=u_{*}+c_{0}-c_{0}\frac{q_{*}}{q_{1}}\ \ ,\ u_{2}=u_{*}+c_{0}-c_{0}\frac{q_{*}}{q_{2}}\ 
,\end{displaymath}  which gives  
\begin{displaymath} 
\frac{u_{2}-u_{1}}{q_{2}-q_{1}}\ =\ \frac{c_{0}q_{*}}{q_{1}q_{2}}\ .\end{displaymath} 
 We have obtained
\begin{displaymath} 
\left({\frac{u_{2}-u_{1}}{q_{2}-q_{1}}\ =}\right) \ c_{0}\ \frac{1}{\sqrt{q_{1}q_{2}}}\ 
=\ c_{0}\ \frac{1}{q_{0}}\ =\ c_{0}\ \frac{q_{*}}{q_{1}q_{2}}\ ,\ {\mathrm{hence}}\ 
q_{*}=\sqrt{q_{1}q_{2}}=q_{0}\ .\end{displaymath}  Thus $u_{*}=u_{0}$ 
and $(u_{*}=)u_{1}+c_{0}-c_{0}\frac{q_{1}}{q_{*}}=u_{1}-c_{0}+c_{0}\frac{q_{*}}{q_{1}}$, 
which gives $2q_{*}q_{1}-q_{1}^{2}-q_{*}^{2}=0,\ $that is $\left({q_{*}-q_{1}}\right) ^{2}=0,\ $or 
$q_{1}=q_{*}.$ The same arguments give $q_{2}=q_{*}.$ Our hypothetical 
shock is reduced to one point, which means obviously that no shock occurs.\ 
\par
{\hskip 1.8em}The wave velocity $A$ is close to $c_{0},$ but a priori 
different since $A-c_{0}=u_{2}$ for the front part of the wave for instance. 
Depending on the (still realistic) values of the parameters of the model, 
$u_{2}$ can be negative or positive, and a possible constraint on the 
parameters may consist into setting $u_{2}=0.\ $The arguments for this 
approach are settled in the reference [2].\ \par
{\hskip 1.8em}This construction of a wave should also feed the discussion 
between linear or non linear models (see for instance [5]). It seems, 
as a conclusion, that the difference is not so important when a source 
term is present. By taking in account a non linear diffusion term instead 
of a friction term, other travelling waves may be expected. The main 
 conclusion is that the linearity is naturally induced by the shape in 
a wind instrument.\ \par
{\hskip 1.8em}This paper is also available on the website "conservation 
laws preprint server", 2008-010.\ \par
\ \par
{\hskip 18.0em}References\ \par
[1] R.F.Dressler, {\it{Mathematical solution of the problem of Roll waves 
in Inclined Open Channel,        }}       Comm. Pure Appl. Math. Vol.II, 
N{\char'27}2-3, pp149-194 (1949).\ \par
[2] J.Kergomard, J.D.Polack, J.Gilbert,{\it{Vitesse de propagation d'une 
onde plane impulsionelle dans un       tuyau sonore,}} J.Acoustique 4 
pp467-483 (1991).\ \par
[3] A.Y.LeRoux, M.N.LeRoux,{\it{Source waves}}, On the Conservation Laws 
Preprint Server\ \par
{\hskip 1.8em}http:$/$www.math.ntnu.no$/$conservation$/$2004$/$045.html 
\ \par
[4] A.Y.LeRoux, M.N.LeRoux, {\it{A mathematical model for Rogue Waves, 
using SaintVenant equations with friction}},  On the Conservation Laws 
Preprint Server \ \par
{\hskip 1.8em}http:$/$www.math.ntnu.no$/$conservation$/$2005$/$048.html\ 
\par
[5] S.Scheichl, {\it{Linear and nonlinear propagation of higher order 
modes in hard-walled circular             ducts containing a real gas,}} 
by  J.Acous.Soc.Am. 117(4) pp1806-1827 (2005).\ \par
\ \par

\end{document}